\documentclass[graybox]{svmult}

\usepackage{type1cm}        
\usepackage{makeidx}         
\usepackage{graphicx}        
\usepackage{multicol}        
\usepackage[bottom]{footmisc}

\usepackage{newtxtext}       %
\usepackage[varvw]{newtxmath}       

\makeindex             

\begin{document}

\title*{Discrete octonionic analysis: a unified approach to the split-octonionic and classical settings}
\titlerunning{Discrete octonionic analysis: split-octonions and classical settings}
\author{Rolf S\"oren Krau\ss har\orcidID{0000-0002-0915-1898} and\\ Anastasiia Legatiuk\orcidID{0000-0003-3903-6281} and 
\\ Dmitrii Legatiuk\orcidID{0000-0002-0028-5793}}
\institute{Rolf S\"oren Krau\ss har \at University of Erfurt, Nordh\"auser Str. 63, 99089 Erfurt, Germany  \email{soeren.krausshar@uni-erfurt.de}
\and Anastasiia Legatiuk \at  University of Erfurt, Nordh\"auser Str. 63, 99089 Erfurt, Germany \email{anastasiia.legatiuk@uni-erfurt.de}
\and Dmitrii Legatiuk \at  University of Erfurt, Nordh\"auser Str. 63, 99089 Erfurt, Germany \email{dmitrii.legatiuk@uni-erfurt.de}}
%
%
\maketitle

\abstract{Various problems of mathematical physics consider octonions and split-octonions as a mathematical structure, which underpins the eight-dimensional nature of these problems. Therefore, it is not surprising that octonionic analysis has become an area of active research in recent years. One of the main goals of octonionic analysis is to develop tools of an octonionic operator calculus for solving boundary value problems of mathematical physics that benefit from the use of the octonionic structure. However, when we want to apply the operator calculus in practice, it becomes evident that adequate discrete counterparts of continuous constructions need to be defined. In previous works, we have proposed several approaches to discretise the classical continuous octonionic analysis. However, the split-octonionic case, which is particularly important for practical applications concretely investigated in the last years, has not been considered until now. Therefore, one of the goals of this paper is to explain how to particularly address the discrete split-octonionic setting. Additionally, we propose a general umbrella to cover all different discrete octonionic settings in one unified approach that also encompasses the different eight-dimensional algebraic structures.}

\section{Introduction}

Recent developments in mathematical physics, such as for instance the modelling of dyonic plasma dynamics, indicated that it is advantageous to treat the related dyonic Maxwell equation and the more sophisticated dyonic magneto-hydrodynamic equations using the structure of eight-dimensional non-associative real algebras. The typical choice for such non-associative real algebra is for example the Cayley octonions or, alternatively, the so-called algebra of octons which represents an $8$-dimensional non-associative subalgebra of the $16$-dimensional sedenions, see for example \cite{KTD2021,DZ,KTD2020}. However, in particular, the use of the hyperbolic octonions which are often also called split-octonions offer an elegant way to formulate and to treat the fluid Maxwell equations, see \cite{BS2016,DT2016,GG2022} for details. As explicitly shown in \cite{DT2016}, the relations between the fluid field equations can be established by the simple application of the multiplication operation in the split-octonions. This result clearly shows that the embedding of this problem into the algebraic structure of the split-octonions is a canonical and therefore adequate choice. Furthermore, split-octonions additionally allow us to elegantly combine the field equations of electromagnetism with linear gravity, see for example \cite{DTT2013}.\par
In recent years, a lot of effort has been put in the fundamental development of a function theory and operator theory in the classical octonionic setting, see \cite{KFVR2024} and the abundant list of references therein. In contrast, the case of split-octonions has not received a sufficient attention until now. Evidently, advancing the results related to split-octonions would enhance the analysis of the applied problems mentioned above, as well as the development of solution methods for such problems.\par
In works \cite{Krausshar_1,Krausshar_2,Krausshar_4,Krausshar_5}, ideas towards discretisation of octonionic differential and integral operators based on several discretisation approaches used previously also in the associative Clifford algebra setting have been presented. In particular, a discrete octononic analysis based on the forward, backward and central discrete Cauchy-Riemann operators has been proposed in \cite{Krausshar_2}. A Weyl calculus perspective on the discrete octonionic analysis in bounded domains in $\mathbb{R}^{8}$ and for half-spaces has been presented in \cite{Krausshar_4}. This approach adapts ideas from the associated setting of discrete Clifford analysis developed in \cite{Brackx_1, Cerejeiras}. For the ideas on using forward and backward discretisations to address bounded domains we refer to \cite{RZ}.\par
In this paper extends the ideas of discrete octonionic analysis to the case of split-octonions. However, to avoid the repetition of standard proofs, we propose a general umbrella that covers on the one hand all different discretisations mentioned above, and on the other hand also all the different octonionic algebras. Consequently, this paper aims to introduce this general umbrella, which enables a unified treatment of the split-octonions and the classical octonions in different discretisation settings. The advantage then consists of the applicability of this general toolkit to all the different kinds of eight-dimensional problems. In particular, the problems related to the dyonic plasma equations and dyonic magneto-hydrodynamic equations could also be addressed.\par      
The paper is organised as follows: Section~\ref{Section:preliminaries} provides the basic notions and definitions of the classical octonions and the split-octonions together with definitions of discrete Cauchy-Riemann operators and their fundamental solutions. After that, a unified approach covering all operators in all these different eight-dimensional real algebras is introduced in Section~\ref{Section:unified}. The key part of this unified approach is the definition of two different abstract associator expressions executed over index sets which take into account the particular multiplication rules in the algebra. This allows us to define a general abstract Stokes operator together with abstract boundary operators. In Section~\ref{Section:discrete_analysis}, we then explain how these are used to set up general Stokes and Borel Pompeiu formulae (including generalisations of Cauchy's formula) covering all these algebraic settings and different discretisations of the associated Cauchy-Riemann operators in one. Finally, we address the particular half-space case in the octonions and split-octonions.\par

\section{Preliminaries}\label{Section:preliminaries}

In this section, we briefly summarise some basic ideas of discrete octonionic analysis. For the details we refer the reader to our previous works  \cite{Krausshar_1,Krausshar_2,Krausshar_4,Krausshar_5}. As commonly used, a vector $\mathbf{x} = (x_{0}, x_{1},\ldots, x_{7})$ from $\mathbb{R}^8$ constituted by the real components $x_0, x_1, \ldots, x_7$ can be associated with a classical octonionic number or split-octonionic number simply as follows 
\begin{equation*}
x = x_{0}\mathbf{e}_{0}+x_{1}\mathbf{e}_{1}+x_{2}\mathbf{e}_{2}+x_{3}\mathbf{e}_{3}+x_{4}\mathbf{e}_{4}+x_{5}\mathbf{e}_{5}+x_{6}\mathbf{e}_{6}+x_{7}\mathbf{e}_{7}. 
\end{equation*}
The multiplication rules for the basis elements $\mathbf{e}_{i}$ are provided in Table~\ref{Table:octonions} for the case of classical octonions (left) and for the case of split-octonions (right), which are also considered in this paper.  The main difference between octonions $\mathbb{O}$ and split-octonions $\mathbb{O}'$ is that split-octonions contain zero divisors, which require some additional attention in the continuous setting.\par
\begin{table}
\caption{Multiplication table for real octonions $\mathbb{O}$ (left) and for split-octonions $\mathbb{O}'$ (right)}
\label{Table:octonions}
\begin{center}
	\begin{tabular}{|c|rrrrrrr|}
\hline
$\cdot$ & $\mathbf{e}_{1}$&  $\mathbf{e}_{2}$ & $\mathbf{e}_{3}$ & $\mathbf{e}_{4}$ & $\mathbf{e}_{5}$ & $\mathbf{e}_{6}$  & $\mathbf{e}_{7}$ \\ \hline
$\mathbf{e}_{1}$  &  $-1$ &  $\mathbf{e}_{4}$ & $\mathbf{e}_{5}$ & $-\mathbf{e}_{2}$ &$-\mathbf{e}_{3}$ & $-\mathbf{e}_{7}$ & $\mathbf{e}_{6}$ \\
$\mathbf{e}_{2}$ &  $-\mathbf{e}_{4}$&   $-1$ & $\mathbf{e}_{6}$ & $\mathbf{e}_{1}$ & $\mathbf{e}_{7}$ & $-\mathbf{e}_{3}$ & $-\mathbf{e}_{5}$ \\
$\mathbf{e}_{3}$ &  $-\mathbf{e}_{5}$& $-\mathbf{e}_{6}$ & $-1$  & $-\mathbf{e}_{7}$ & $\mathbf{e}_{1}$  & $\mathbf{e}_{2}$  & $\mathbf{e}_{4}$ \\
$\mathbf{e}_{4}$ &  $\mathbf{e}_{2}$ & $-\mathbf{e}_{1}$ & $\mathbf{e}_{7}$ & $-1$  &$-\mathbf{e}_{6}$ & $\mathbf{e}_{5}$  & $-\mathbf{e}_{3}$\\
$\mathbf{e}_{5}$ &  $\mathbf{e}_{3}$ & $-\mathbf{e}_{7}$ & $-\mathbf{e}_{1}$&  $\mathbf{e}_{6}$&  $-1$ & $-\mathbf{e}_{4}$ & $\mathbf{e}_{2}$ \\
$\mathbf{e}_{6}$ &  $\mathbf{e}_{7}$ &  $\mathbf{e}_{3}$ & $-\mathbf{e}_{2}$& $-\mathbf{e}_{5}$ & $\mathbf{e}_{4}$ & $-1$   & $-\mathbf{e}_{1}$ \\
$\mathbf{e}_{7}$ & $-\mathbf{e}_{6}$ &  $\mathbf{e}_{5}$ & $-\mathbf{e}_{4}$ & $\mathbf{e}_{3}$ & $-\mathbf{e}_{2}$ & $\mathbf{e}_{1}$  & $-1$ \\ \hline 	
\end{tabular}
 \qquad 
	\begin{tabular}{|c|rrrrrrr|}
\hline
$\cdot$ & $\mathbf{e}_{1}$&  $\mathbf{e}_{2}$ & $\mathbf{e}_{3}$ & $\mathbf{e}_{4}$ & $\mathbf{e}_{5}$ & $\mathbf{e}_{6}$  & $\mathbf{e}_{7}$ \\ 
\hline
$\mathbf{e}_{1}$  &  $-1$ &  $\mathbf{e}_{3}$ & $\mathbf{e}_{2}$ & $-\mathbf{e}_{5}$ & $\mathbf{e}_{4}$ & $-\mathbf{e}_{7}$ & $\mathbf{e}_{6}$ \\
$\mathbf{e}_{2}$ &  $-\mathbf{e}_{3}$&   $-1$ & $\mathbf{e}_{1}$ & $-\mathbf{e}_{6}$ & $\mathbf{e}_{7}$ & $\mathbf{e}_{4}$ & $-\mathbf{e}_{5}$ \\
$\mathbf{e}_{3}$ &  $\mathbf{e}_{2}$ & $-\mathbf{e}_{1}$ & $-1$  & $-\mathbf{e}_{7}$ & $\mathbf{e}_{6}$  & $\mathbf{e}_{5}$ & $\mathbf{e}_{4}$ \\
$\mathbf{e}_{4}$ &  $\mathbf{e}_{5}$ & $\mathbf{e}_{6}$ & $\mathbf{e}_{7}$ & $1$  & $\mathbf{e}_{1}$ & $\mathbf{e}_{2}$  & $-\mathbf{e}_{3}$\\
$\mathbf{e}_{5}$ &  $-\mathbf{e}_{4}$ & $-\mathbf{e}_{7}$ & $\mathbf{e}_{6}$ &  $-\mathbf{e}_{1}$&  $1$ & $\mathbf{e}_{3}$ & $-\mathbf{e}_{2}$ \\
$\mathbf{e}_{6}$ &  $\mathbf{e}_{7}$ &  $-\mathbf{e}_{4}$ & $-\mathbf{e}_{5}$ & $-\mathbf{e}_{2}$ & $-\mathbf{e}_{3}$ & $1$ & $\mathbf{e}_{1}$ \\
$\mathbf{e}_{7}$ & $-\mathbf{e}_{6}$ &  $\mathbf{e}_{5}$ & $-\mathbf{e}_{4}$ & $-\mathbf{e}_{3}$ & $\mathbf{e}_{2}$ & $-\mathbf{e}_{1}$  & $1$ \\ \hline 	
\end{tabular}
\end{center}
\end{table}\par
As it can clearly be seen from Table~\ref{Table:octonions}, the multiplication of octonions and split-octonions is closed. However, in both cases it is non-associative, because we have the relations
\begin{equation*}
\left\{\begin{array}{rcll}
(\mathbf{e}_{i}\mathbf{e}_{j})\mathbf{e}_{k} & = & -\mathbf{e}_{i}(\mathbf{e}_{j}\mathbf{e}_{k}) & \mbox{for mutually distinct, non-zero } i,j,k \mbox{ and } \mathbf{e}_{i}\mathbf{e}_{j} \neq \pm \mathbf{e}_{k},  \\[2mm]
(\mathbf{e}_{i}\mathbf{e}_{j})\mathbf{e}_{k} & = & \mathbf{e}_{i}(\mathbf{e}_{j}\mathbf{e}_{k}) & \mbox{otherwise}.
\end{array}\right.
\end{equation*}\par
Next, we introduce the unbounded uniform lattice $h\mathbb{Z}^{8}$ with $h>0$ in the usual way
\begin{equation*}
h \mathbb{Z}^{8} :=\left\{\mathbf{x} \in {\mathbb{R}}^{8}\,|\, \mathbf{x} = (m_{0}h, m_{1}h,\ldots, m_{7}h), m_{j} \in \mathbb{Z}, j=0,1,\ldots,7\right\}.
\end{equation*}
Additionally, we introduce the upper and lower half-spaces:
\begin{equation*}
\begin{array}{rcl}
h\mathbb{Z}_{+}^{8} & := & \left\{(h\underline{m},hm_{7})\colon \underline{m}\in\mathbb{Z}^{7},m_{7}\in\mathbb{Z}_{+}\right\}, \\
h\mathbb{Z}_{-}^{8} & := & \left\{(h\underline{m},hm_{7})\colon \underline{m}\in\mathbb{Z}^{7},m_{7}\in\mathbb{Z}_{-}\right\}.
\end{array}
\end{equation*}
Furthermore, for shortening the notations, the arguments of functions will be written as $mh$ instead of the full version $m_{0}h, m_{1}h,\ldots, m_{7}h$.\par
Next, the classical forward and backward differences $\partial_{h}^{\pm j}$ are given by
\begin{equation}
\label{Finite_differences}
\begin{array}{lcl}
\partial_{h}^{+j}f(mh) & := & h^{-1}(f(mh+\mathbf{e}_jh)-f(mh)), \\
\partial_{h}^{-j}f(mh) & := & h^{-1}(f(mh)-f(mh-\mathbf{e}_jh)),
\end{array}
\end{equation}
for discrete functions $f(mh)$ with $mh\in h\mathbb{Z}^{8}$. In this paper, we consider functions defined on $\Omega_{h} \subset  h\mathbb{Z}^{8}$ and taking values in octonions $\mathbb{O}$ or split-octonions $\mathbb{O}'$, respectively. We will always indicate clearly which setting is considered at the moment and, moreover, we will also underline the differences, if any, appearing when considering $\mathbb{O}'$ instead of $\mathbb{O}$.\par
By the help of finite difference operators~(\ref{Finite_differences}), we can now introduce a {\itshape discrete octonionic forward Cauchy-Riemann operator} $D^{+}\colon l^{p}(\Omega_{h},\mathbb{O})\to l^{p}(\Omega_{h},\mathbb{O})$ and a {\itshape discrete octonionic backward Cauchy-Riemann operator} $D^{-}\colon l^{p}(\Omega_{h},\mathbb{O})\to l^{p}(\Omega_{h},\mathbb{O})$ as follows
\begin{equation}
\label{Cauchy_Riemann_operators_discrete}
D^{+}_{h}:=\sum_{j=0}^{7} \mathbf{e}_j\partial_{h}^{+j}, \quad D^{-}_{h}:=\sum_{j=0}^{7} \mathbf{e}_j\partial_{h}^{-j},
\end{equation}
for $1\leq p<\infty$. These discrete octonionic forward and backward Cauchy-Riemann operators give rise to {\itshape discrete octonionic forward and backward monogenic functions}:
\begin{definition}\label{Definition:discrete_monogenic}
A function $f\in l^{p}(\Omega_{h},\mathbb{O})$ is called {\itshape discrete octonionic left forward monogenic} if $D_{h}^{+}f=0$ in $\Omega_{h}$. Respectively, a function $f\in l^{p}(\Omega_{h},\mathbb{O})$ is called {\itshape discrete octonionic left backward monogenic} if $D_{h}^{-}f=0$ in $\Omega_{h}$.
\end{definition}
Discrete Cauchy-Riemann operators and discrete monogenic functions for split-octonions $\mathbb{O}'$ could be defined analogously. No principal differences between both settings are being manifested at this level.\par
For the construction of discrete counterparts of the continuous Borel-Pompeiu and Cauchy formulae, it is necessary to work with the discrete fundamental solution of discrete octonionic Cauchy-Riemann operators. We recall its definition for the case of octonions, while the split-octonionic case is completely analogous:
\begin{definition}
The function $E_{h}^{+}\colon h\mathbb{Z}^{8} \rightarrow \mathbb{O}$ is called a discrete fundamental solution of $D_{h}^{+}$ if it satisfies
\begin{equation*}
D_{h}^{+}E_{h}^{+} =\delta_h = \begin{cases} h^{-8}, & \mbox{for } mh=0,\\ 
0,& \mbox{for } mh\neq 0,
\end{cases}
\end{equation*} 
for all grid points $mh$ of $h\mathbb{Z}^{8}$. Analogously, the function $E_{h}^{-}\colon h\mathbb{Z}^{8} \rightarrow \mathbb{O}$ is called a discrete fundamental solution of $D_{h}^{-}$ if it satisfies
\begin{equation*}
D_{h}^{-}E_{h}^{-} =\delta_h = \begin{cases} h^{-8}, & \mbox{for } mh=0,\\ 
0,& \mbox{for } mh\neq 0,
\end{cases}
\end{equation*} 
for all grid points $mh$ of $h\mathbb{Z}^{8}$.
\end{definition}\par
The construction of the discrete fundamental solutions $E_{h}^{\pm}$ follows the classical scheme of using the discrete Fourier transform, see for example \cite{Stummel}. The octonionic case has been considered in \cite{Krausshar_2}, and the discrete fundamental solution $E_{h}^{\pm}$ is expressed as follows
\begin{equation}
\label{Discrete_fundamental_solution}
E_{h}^{\pm}=\sum\limits_{j=0}^{7}\mathbf{e}_{j}\mathcal{R}_{h} \mathcal{F} \left( \frac{\xi_{h}^{\pm j}}{d^2} \right),
\end{equation}
where $\xi_{h}^{\pm j}=\mp h^{-1}\left( 1-e^{\mp ih\xi_j}\right)$ are the Fourier symbols of the finite difference operators~(\ref{Finite_differences}), $d^2=\frac{4}{h^2}\sum\limits_{j=0}^{7} \sin^2\left(\frac{\xi_j h}{2}\right)$ is the symbol of the classical star-Laplacian, and $\mathcal{R}_{h} \mathcal{F}$ is the restriction of the classical continuous Fourier transform $\mathcal{F}$ to the lattice. Again, the split-octonionic structure does not affect the construction principle of $E_{h}^{\pm}$. In fact, the particular effect of the split-octonionic structure will be visible in the multiplication rules for $\mathbf{e}_{j}$ appearing in~(\ref{Discrete_fundamental_solution}), but this will be immediately clear from the context. Therefore, we will use the discrete fundamental solutions $E_{h}^{\pm}$ for $\mathbb{O}$ and $\mathbb{O}'$ without introducing a further particular notations.\par

\section{A unified approach to the discrete setting}\label{Section:unified}

In previous papers we have developed a discrete octonionic analysis in two general settings: one is based on the use of discrete octonionic forward and backward Cauchy-Riemann operators, and another one is based on the Weyl calculus approach adapted from the discrete Clifford analysis, see for example \cite{Brackx_1,CKKS,Cerejeiras,Cerejeiras_2,Vaz}. Evidently, split-octonions could also be addressed in both settings. However, as we have already mentioned and as we will observe more deeply in the next section, the split-octonionic setting neither changes critically the structure of proof nor the structure of the arising final formulae. Therefore, our goal in this short section is to propose a unified approach to generally present  the results in the discrete octonionic setting.\par
A distinct feature of octonionic analysis is the presence of so called {\itshape associator}, which in the continuous case defined as $[a,b,c]:=(ab)c-a(bc)$ for octonions $a,b,c$ \cite{XL2002}. In the discrete setting, the associator has been overlooked during the first attempts to construct a consistent theory \cite{Krausshar_1,Krausshar_2,Krausshar_3}. However, it has been recovered in \cite{Krausshar_5} for the case of discrete octonionic forward and backward Cauchy-Riemann operators. In \cite{Krausshar_4} it has also been recovered for the Weyl calculus-based approach. The key idea of the Weyl calculus-based approach is the splitting of basis elements $\mathbf{e}_{k}$, $k=0,1,\ldots,7$, into positive and negative directions $\mathbf{e}_{k}^+$ and $\mathbf{e}_{k}^-$, $k=0,1,\ldots,7$, i.e., $\mathbf{e}_{k}=\mathbf{e}_{k}^{+}+\mathbf{e}_{k}^-$. This splitting must satisfy the following relations
\begin{equation*}
\left\{
\begin{array}{ccc}
\mathbf{e}_j^{-}\mathbf{e}_k^{-}+\mathbf{e}_k^{-}\mathbf{e}_j^{-} &=&0, \\
\mathbf{e}_j^{+}\mathbf{e}_k^{+}+\mathbf{e}_k^{+}\mathbf{e}_j^{+}&=&0, \\
\mathbf{e}_j^{+}\mathbf{e}_k^{-} + \mathbf{e}_k^{-}\mathbf{e}_j^{+}&=&-\delta_{jk},
\end{array}
\right.
\end{equation*}
where $\delta_{jk}$ is the Kronecker delta symbol. In the discrete setting, this splitting give rise to the pair of discrete Cauchy-Riemann operators
\begin{equation*}
D^{+-}_{h}:=\sum_{j=0}^{7} \mathbf{e}_{j}^{+}\partial_{h}^{+j}+\mathbf{e}_{j}^{-}\partial_{h}^{-j}, \quad D^{-+}_{h}:=\sum_{j=0}^{7} \mathbf{e}_{j}^{+}\partial_{h}^{-j}+\mathbf{e}_{j}^{-}\partial_{h}^{+j},
\end{equation*}
which factorise the classical star-Laplacian, see \cite{CKKS,Cerejeiras} for details.\par
To cover the cases of forward and backward Cauchy-Riemann operators and the Weyl calculus-based approach, let us introduce the following {\itshape abstract discrete associators}
\begin{equation}
\label{Abstrac_associators}
\begin{array}{ccl}
\mathcal{A}_{1}(f,g,\mathcal{N},\mathcal{K}) & := & \displaystyle 2\sum_{m\in \mathcal{N}} \sum\limits_{s=1}^{7}  \sum_{\stackrel{i=1}{i\in I_{s}}}^{7} \sum_{\stackrel{j=1,j\in \mathcal{K}_{s}}{j\neq i}}^{7} \sum_{\stackrel{k=1}{k\notin \mathcal{K}_{s}}}^{7} \left[g_{i}(mh)\mathbf{e}_{i}\left(\partial_{h}^{-j}\mathbf{e}_{j}f_{k}(mh)\mathbf{e}_{k}\right)\right]h^{8}, \\
\mathcal{A}_{2}(f,g,\mathcal{N},\mathcal{K}) & := & \displaystyle 2\sum_{m\in \mathcal{N}} \sum\limits_{s=1}^{7} \sum_{\stackrel{i=1}{i\in I_{s}}}^{7} \sum_{\stackrel{j=1,j\in \mathcal{K}_{s}}{j\neq i}}^{7} \sum_{\stackrel{k=1}{k\notin \mathcal{K}_{s}}}^{7} \left[g_{i}(mh)\mathbf{e}_{i}\left(\partial_{h}^{+j}\mathbf{e}_{j}^{+}f_{k}(mh)\mathbf{e}_{k}\right) \right. \\
& & \displaystyle \left. +g_{i}(mh)\mathbf{e}_{i}\left(\partial_{h}^{-j}\mathbf{e}_{j}^{-}f_{k}(mh)\mathbf{e}_{k}\right)\right]h^{8},
\end{array}
\end{equation}
where $\mathcal{N}$ is the index set of a domain over which the associator should be calculated, e.g. $\mathbb{Z}^{8}$ for the case of the whole space, and $\mathcal{K}=\left\{I,J\right\}$ is the set containing sets of indices for the case of octonions ($I$) and split-octonions ($J$). Precisely speaking, these sets are defined as follows:
\begin{equation*}
I=\left\{I_{1},I_{2},I_{3},I_{4},I_{5},I_{6},I_{7}\right\} \mbox{ and } J=\left\{J_{1},J_{2},J_{3},J_{4},J_{5},J_{6},J_{7}\right\}
\end{equation*}
with
\begin{equation*}
\begin{array}{cclcclcclccl}
I_{1} & := & \left\{1,2,4\right\}, & I_{2} & := & \left\{1,3,5\right\}, & I_{3} & := & \left\{1,6,7\right\}, & I_{4} & := & \left\{2,3,6\right\}, \\
I_{5} & := & \left\{2,5,7\right\}, & I_{6} & := & \left\{3,4,7\right\}, & I_{7} & := & \left\{4,5,6\right\}.
\end{array}
\end{equation*}
and 
\begin{equation*}
\begin{array}{cclcclcclccl}
J_{1} & := & \left\{1,2,3\right\}, & J_{2} & := & \left\{1,4,5\right\}, & J_{3} & := & \left\{1,6,7\right\}, & J_{4} & := & \left\{2,4,6\right\}, \\
J_{5} & := & \left\{2,5,7\right\}, & J_{6} & := & \left\{3,4,7\right\}, & J_{7} & := & \left\{4,5,6\right\}.
\end{array}
\end{equation*}
The abstract discrete associator $\mathcal{A}_{1}$ gives rise to the discrete octonionic forward and backward Cauchy-Riemann operators, while the abstract associator $\mathcal{A}_{2}$ gives rise to the Weyl calculus-based approach. Note that both associators encompass the classical octonionic setting and the split-octonionic setting as particular settings in one generalised formal framework.\par
\begin{remark}
It is necessary to comment on the origin of sets $I$ and $J$. During the proof of the  discrete Stokes' formula in \cite{Krausshar_5}, it has been shown that the non-associative sub-part of octonionic discrete Stokes' formula cannot be written in a compact way with the aid of discrete octonionic forward and backward Cauchy-Riemann operators. Instead, it turned out to be possible to write them more compactly by the help of the index sets $I_{k}$, $k=1,\ldots,7$. Analogously, the split-octonionic case then leads to the index sets $J_{k}$, $k=1,\ldots,7$.
\end{remark}\par
Next, we introduce {\itshape an abstract Stokes' operator} as follows
\begin{equation}
\label{Abstract_Stokes}
\mathfrak{S}(f,g,\mathcal{D}_{1},\mathcal{D}_{2},\mathcal{N}) := \sum_{m\in \mathcal{N}}  \left[ \left( g(mh)\mathcal{D}_{1}\right) f(mh) + g(mh) \left( \mathcal{D}_{2}f(mh) \right) \right] h^8,
\end{equation}
where $\mathcal{N}$ is the index set of a domain over which the formula should be calculated, and $\mathcal{D}_{1},\mathcal{D}_{2}$ is a pair of discrete Cauchy-Riemann operators used to build the discrete octonionic function theory, e.g. $D_{h}^{+}$ and $D_{h}^{-}$ in the current paper or $D_{h}^{-+}$ and $D_{h}^{+-}$ for the Weyl calculus-based approach.\par
By the help of the abstract discrete associators~(\ref{Abstrac_associators}) and the abstract Stokes' operator~(\ref{Abstract_Stokes}), we can now formally write the general discrete Stokes' formula:
\begin{equation*}
\mathfrak{S}(f,g,\mathcal{D}_{1},\mathcal{D}_{2},\mathcal{N}) = \mathcal{A}_{i}(f,g,\mathcal{N},\mathcal{K}), \mbox{ with } i=1,2.
\end{equation*}
This formula now covers simultaneously all the cases we want to consider, namely the Weyl calculus-based approach and the approach based on forward and backward Cauchy-Riemann operators both for octonions $\mathbb{O}$ and split-octonions $\mathbb{O}'$. In the sequel  of the paper, we will follow this unified approach to the discrete setting, because it provides a compact way to present building blocks of a discrete octonionic function theory. Nonetheless, to improve the readability of this paper, we will also present some discrete formulae explicitly.\par 
Finally, we introduce the following {\itshape abstract boundary operators}, which are needed to define discrete formulae over domains with boundary:
\begin{equation}
\label{Abstrac_boundary_operators}
\begin{array}{ccl}
\mathfrak{B}_{1}(f,g,\mathcal{N},L_{1},L_{2}) & := & \displaystyle 2\sum\limits_{\underline{m}\in\mathcal{N}}\left[\sum\limits_{i=1,6}\sum_{\stackrel{k=1}{k\neq i}}^{6}g_{i}(\underline{m}h,L_{1})\mathbf{e}_{i}\left(\mathbf{e}_{7}f_{k}(\underline{m}h,L_{2})\mathbf{e}_{k}\right) \right. \\
& & \displaystyle + \sum\limits_{i=1,6}\sum_{\stackrel{k=1}{k\neq i}}^{6}g_{i}(\underline{m}h,L_{1})\mathbf{e}_{i}\left(\mathbf{e}_{7}f_{k}(\underline{m}h,L_{2})\mathbf{e}_{k}\right) \\
& & \displaystyle \left. + \sum\limits_{i=1,6}\sum_{\stackrel{k=1}{k\neq i}}^{6}g_{i}(\underline{m}h,L_{1})\mathbf{e}_{i}\left(\mathbf{e}_{7}f_{k}(\underline{m}h,L_{2})\mathbf{e}_{k}\right) \right]h^{7}, \\
\mathfrak{B}_{2}(f,g,\mathcal{N},L_{1},L_{2}) & := & \displaystyle 2\sum\limits_{\underline{m}\in\mathcal{N}}\left[\sum_{i=1,6} \sum_{\stackrel{k=1}{k\neq i}}^{6} \left(g_{i}(\underline{m}h,L_{1})\mathbf{e}_{i}\left(\mathbf{e}_{7}^{+}f_{k}(\underline{m}h,L_{2})\mathbf{e}_{k}\right) \right. \right. \\
& & \displaystyle \left. + g_{i}(\underline{m}h,L_{2})\mathbf{e}_{i}\left(\mathbf{e}_{7}^{-}f_{k}(\underline{m}h,L_{1})\mathbf{e}_{k}\right)\right) \\
& & \displaystyle \sum_{i=2,5} \sum_{\stackrel{k=1}{k\neq i}}^{6} \left(g_{i}(\underline{m}h,L_{1})\mathbf{e}_{i}\left(\mathbf{e}_{7}^{+}f_{k}(\underline{m}h,L_{2})\mathbf{e}_{k}\right) \right. \\
& & \displaystyle \left. + g_{i}(\underline{m}h,L_{2})\mathbf{e}_{i}\left(\mathbf{e}_{7}^{-}f_{k}(\underline{m}h,L_{1})\mathbf{e}_{k}\right)\right)\\
& & \displaystyle \sum_{i=3,4} \sum_{\stackrel{k=1}{k\neq i}}^{6} \left(g_{i}(\underline{m}h,L_{1})\mathbf{e}_{i}\left(\mathbf{e}_{7}^{+}f_{k}(\underline{m}h,L_{2})\mathbf{e}_{k}\right) \right. \\
& & \displaystyle \left. \left. + g_{i}(\underline{m}h,L_{2})\mathbf{e}_{i}\left(\mathbf{e}_{7}^{-}f_{k}(\underline{m}h,L_{1})\mathbf{e}_{k}\right)\right)\right] h^{7},
\end{array}
\end{equation}
where $\mathcal{N}$ is the index set of a domain over which the boundary operator should be calculated, $L_{1}$ and $L_{2}$ are the indices of boundary layers. In the discrete setting, the boundary of a domain is represented by a three-layer structure with the middle layer being a \textquotedblleft real\textquotedblright\, boundary of the domain, we refer to \cite{Cerejeiras,Krausshar_4} for a complete discussion on the discrete boundary. Again $\mathfrak{B}_{1}$ corresponds to the case of discrete octonionic forward and backward Cauchy-Riemann operators and $\mathfrak{B}_{2}$ to the case of the  Weyl calculus-based approach.\par

\section{Discrete octonionic analysis: classical and split-octonionic settings}\label{Section:discrete_analysis}

In \cite{Krausshar_1,Krausshar_2} we have developed a discrete counterpart of octonionic analysis including Stokes', Borel-Pompeiu, and Cauchy formulae. However, during this development an important algebraic property of octonions has been overseen. This feature  has been taken into account in \cite{Krausshar_5}, where the discrete Stokes' formula for the whole space has been revisited and corrected. Similarly, in \cite{Krausshar_4} we revisited  the discrete octonionic analysis for bounded domains in the framework of the  Weyl calculus approach. These considerations exposed the necessity to also revisit the results discussed in \cite{Krausshar_2} that were related to the discrete octonionic analysis for half-spaces $h\mathbb{Z}_{+}^{8}$ and $h\mathbb{Z}_{-}^{8}$. Moreover, we will use the unified approach proposed in the previous section to treat simultaneously the octonionic and split-octonionic case. After that, some formulae will be presented in explicit form to exhibit clearly what needs to be calculated when using the discrete operator calculus in practice.\par
We start by presenting the discrete Stokes' formula for the whole space $h\mathbb{Z}^{8}$ in an abstract way:
\begin{theorem}
The discrete Stokes' formula for the whole lattice $h\mathbb{Z}^{8}$ is given by
\begin{equation*}
\mathfrak{S}(f,g,D_h^{+},D_h^{-},\mathbb{Z}^{8}) = \mathcal{A}_{1}(f,g,\mathbb{Z}^{8},I)
\end{equation*}
for the classical octonions $\mathbb{O}$, and
\begin{equation*}
\mathfrak{S}(f,g,D_h^{+},D_h^{-},\mathbb{Z}^{8}) = \mathcal{A}_{1}(f,g,\mathbb{Z}^{8},J)
\end{equation*}
for the split-octonions $\mathbb{O}'$. These formulae hold for all discrete functions $f$ and $g$ such that the series converge.
\end{theorem}
It is important to remark that the structure of the discrete Stokes' formulae presented in this theorem makes clear, that the split-octonionic structure only affects the sets of indices appearing in the associator. Hence, working with the unified approach is beneficial here for avoiding repeated analogous proofs and calculations.\par
A next step is to extend this result to the case of half-spaces $h\mathbb{Z}_{+}^{8}$ and $h\mathbb{Z}_{-}^{8}$. The critical difference to the case of the whole space is the presence of boundary terms on three boundary layers, see \cite{Cerejeiras} for a detailed discussion:
\begin{itemize}
\item $m_{7}=0$ and $m_{7}=1$ for the upper half-space $h\mathbb{Z}_{+}^{8}$, and
\item $m_{7}=0$ and $m_{7}=-1$ for the lower half-space $h\mathbb{Z}_{-}^{8}$.
\end{itemize}
By using this information and the abstract boundary operator $\mathfrak{B}_{1}$ from~(\ref{Abstrac_boundary_operators}), we can now introduce the discrete octonionic and split-octonionic Stokes' formulae for the half-spaces:
\begin{theorem}
The discrete Stokes' formula for the upper half-space $h\mathbb{Z}_{+}^{8}$ is given by
\begin{equation*}
\mathfrak{S}(f,g,D_h^{+},D_h^{-},\mathbb{Z}_{+}^{8}) = \mathcal{A}_{1}(f,g,\mathbb{Z}_{+}^{8},I) + \mathfrak{B}_{1}(f,g,\mathbb{Z}^{7},h,0)
\end{equation*}
for the classical octonions $\mathbb{O}$, and
\begin{equation*}
\mathfrak{S}(f,g,D_h^{+},D_h^{-},\mathbb{Z}_{+}^{8}) = \mathcal{A}_{1}(f,g,\mathbb{Z}_{+}^{8},J) + \mathfrak{B}_{1}(f,g,\mathbb{Z}^{7},h,0)
\end{equation*}
for the split-octonions $\mathbb{O}'$. Analogously, the discrete Stokes' formula for the lower half-space $h\mathbb{Z}_{-}^{8}$ is given by
\begin{equation*}
\mathfrak{S}(f,g,D_h^{+},D_h^{-},\mathbb{Z}_{-}^{8}) = \mathcal{A}_{1}(f,g,\mathbb{Z}_{-}^{8},I) + \mathfrak{B}_{1}(f,g,\mathbb{Z}^{7},0,-h)
\end{equation*}
for the classical octonions $\mathbb{O}$, and
\begin{equation*}
\mathfrak{S}(f,g,D_h^{+},D_h^{-},\mathbb{Z}_{-}^{8}) = \mathcal{A}_{1}(f,g,\mathbb{Z}_{-}^{8},J) + \mathfrak{B}_{1}(f,g,\mathbb{Z}^{7},0,-h)
\end{equation*}
for the split-octonions $\mathbb{O}'$. These formulae hold for all discrete functions $f$ and $g$ under the condition that the series converge. 
\end{theorem}
\begin{proof}
The classical octonionic case has been proved in \cite{Krausshar_5}. The split-octonionic case is proved analogously and the only difference appears in the index sets $J_{k}$, $k=1,\ldots,7$, which comes from the different multiplication rules as indicated in Table~\ref{Table:octonions}. After working out carefully these index sets, the final expression for the associator $\mathcal{A}_{1}$ is obtained. Again it should be noted that the split-octonionic structure affects only the index sets in the associator. The remaining steps of the proof from \cite{Krausshar_5} follow immediately.
\end{proof}\par
After introducing the discrete Stokes' formulae for the half-spaces, we can now define the discrete Borel-Pompeiu formulae by substituting the function $g$ by the discrete fundamental solution $E_{h}^{+}(\cdot -mh)$ in the discrete Stokes' formulae. Thus, by the help of the unified approach, we obtain the following theorem:
\begin{theorem}\label{Borel-Pompeiu_abstract}
Let $E_{h}^{+}$ be the discrete fundamental solution to the discrete Cauchy-Riemann operator $D_{h}^{+}$. Then the discrete octonionic Borel-Pompeiu formula for the upper half-space $h\mathbb{Z}_{+}^{8}$ is given by
\begin{equation*}
\mathfrak{S}(f,E_{h}^{+}(nh-mh),D_h^{+},D_h^{-},\mathbb{Z}_{+}^{8}) = \mathcal{A}_{1}(f,g,\mathbb{Z}_{+}^{8},I) + \mathfrak{B}_{1}(f,g,\mathbb{Z}^{7},h,0),
\end{equation*}
and the discrete split-octonionic Borel-Pompeiu formula for the upper half-space $h\mathbb{Z}_{+}^{8}$ is given by
\begin{equation*}
\mathfrak{S}(f,E_{h}^{+}(nh-mh),D_h^{+},D_h^{-},\mathbb{Z}_{+}^{8}) = \mathcal{A}_{1}(f,g,\mathbb{Z}_{+}^{8},J) + \mathfrak{B}_{1}(f,g,\mathbb{Z}^{7},h,0).
\end{equation*}
Analogously, the discrete octonionic Borel-Pompeiu formula for the lower half-space $h\mathbb{Z}_{-}^{8}$ is given by
\begin{equation*}
\mathfrak{S}(f,E_{h}^{+}(nh-mh),D_h^{+},D_h^{-},\mathbb{Z}_{-}^{8}) = \mathcal{A}_{1}(f,g,\mathbb{Z}_{-}^{8},I) + \mathfrak{B}_{1}(f,g,\mathbb{Z}^{7},0,-h),
\end{equation*}
and the discrete split-octonionic Borel-Pompeiu formula for the lower half-space $h\mathbb{Z}_{-}^{8}$ is given by
\begin{equation*}
\mathfrak{S}(f,E_{h}^{+}(nh-mh),D_h^{+},D_h^{-},\mathbb{Z}_{-}^{8}) = \mathcal{A}_{1}(f,g,\mathbb{Z}_{-}^{8},J) + \mathfrak{B}_{1}(f,g,\mathbb{Z}^{7},0,-h).
\end{equation*}
These formulae hold for all discrete functions $f$ and $g$ such that the series converge. 
\end{theorem}
To illustrate the effect of the abstract operator concretely, we now present the discrete Borel-Pompeiu formulae for the upper half-space in explicit form:
\begin{theorem}\label{Borel-Pompeiu_explicit}
Let $E_{h}^{+}$ be the discrete fundamental solution to the discrete Cauchy-Riemann operator $D_{h}^{+}$. Then the discrete octonionic Borel-Pompeiu formula for the upper half-space $h\mathbb{Z}_{+}^{8}$ is given by
\begin{equation*}
\begin{array}{c}
\displaystyle \sum_{n\in \mathbb{Z}_{+}^{8}} E_{h}^{+}(nh-mh) \left[ D_h^{-}f(nh) \right] h^8 -  \\
\displaystyle - 2\sum_{n\in \mathbb{Z}_{+}^{8}} \sum\limits_{s=1}^{7}  \sum_{\stackrel{i=1}{i\in I_{s}}}^{7} \sum_{\stackrel{j=1,j\in I_{s}}{j\neq i}}^{7} \sum_{\stackrel{k=1}{k\notin I_{s}}}^{7} \left[E_{h,i}^{+}(nh-mh)\mathbf{e}_{i}\left(\partial_{h}^{-j}\mathbf{e}_{j}f_{k}(nh)\mathbf{e}_{k}\right)\right]h^{8} \\
\displaystyle + 2\sum\limits_{\underline{n}\in\mathbb{Z}^{7}}\left[\sum\limits_{i=1,6}\sum_{\stackrel{k=1}{k\neq i}}^{6}E_{h,i}^{+}(\underline{n}h-\underline{m}h,h)\mathbf{e}_{i}\left(\mathbf{e}_{7}f_{k}(\underline{n}h,0)\mathbf{e}_{k}\right) \right. \\
\displaystyle + \sum\limits_{i=1,6}\sum_{\stackrel{k=1}{k\neq i}}^{6}E_{h,i}^{+}(\underline{n}h-\underline{m}h,h)\mathbf{e}_{i}\left(\mathbf{e}_{7}f_{k}(\underline{n}h,0)\mathbf{e}_{k}\right) \\
\displaystyle \left. + \sum\limits_{i=1,6}\sum\limits_{\stackrel{k=1}{k\neq i}}^{6}E_{h,i}^{+}(\underline{n}h-\underline{m}h,h)\mathbf{e}_{i}\left(\mathbf{e}_{7}f_{k}(\underline{n}h,0)\mathbf{e}_{k}\right) \right]h^{7} = \left\{\begin{array}{cl}
0, & m\notin \mathbb{Z}_{+}^{8}, \\
-f(mh), & m\in \mathbb{Z}_{+}^{8}.
\end{array} \right.
\end{array}
\end{equation*}
The discrete split-octonionic Borel-Pompeiu formula for the upper half-space $h\mathbb{Z}_{+}^{8}$ is given by
\begin{equation*}
\begin{array}{c}
\displaystyle \sum_{n\in \mathbb{Z}_{+}^{8}} E_{h}^{+}(nh-mh) \left[ D_h^{-}f(nh) \right] h^8 -  \\
\displaystyle - 2\sum_{n\in \mathbb{Z}_{+}^{8}} \sum\limits_{s=1}^{7}  \sum_{\stackrel{i=1}{i\in J_{s}}}^{7} \sum_{\stackrel{j=1,j\in J_{s}}{j\neq i}}^{7} \sum_{\stackrel{k=1}{k\notin J_{s}}}^{7} \left[E_{h,i}^{+}(nh-mh)\mathbf{e}_{i}\left(\partial_{h}^{-j}\mathbf{e}_{j}f_{k}(nh)\mathbf{e}_{k}\right)\right]h^{8} \\
\displaystyle + 2\sum\limits_{\underline{n}\in\mathbb{Z}^{7}}\left[\sum\limits_{i=1,6}\sum_{\stackrel{k=1}{k\neq i}}^{6}E_{h,i}^{+}(\underline{n}h-\underline{m}h,h)\mathbf{e}_{i}\left(\mathbf{e}_{7}f_{k}(\underline{n}h,0)\mathbf{e}_{k}\right) \right. \\
\displaystyle + \sum\limits_{i=1,6}\sum\limits_{\stackrel{k=1}{k\neq i}}^{6}E_{h,i}^{+}(\underline{n}h-\underline{m}h,h)\mathbf{e}_{i}\left(\mathbf{e}_{7}f_{k}(\underline{n}h,0)\mathbf{e}_{k}\right) \\
\displaystyle \left. + \sum\limits_{i=1,6}\sum\limits_{\stackrel{k=1}{k\neq i}}^{6}E_{h,i}^{+}(\underline{n}h-\underline{m}h,h)\mathbf{e}_{i}\left(\mathbf{e}_{7}f_{k}(\underline{n}h,0)\mathbf{e}_{k}\right) \right]h^{7} = \left\{\begin{array}{cl}
0, & m\notin \mathbb{Z}_{+}^{8}, \\
-f(mh), & m\in \mathbb{Z}_{+}^{8}
\end{array} \right.
\end{array}
\end{equation*}
These formulae hold for all discrete functions $f$ and $g$ under the condition that the  the series converge. The discrete Borel-Pompeiu formulae for the lower half-space $h\mathbb{Z}_{-+}^{8}$ are defined analogously.
\end{theorem}\par
Now, comparing the explicit formulae from this theorem with the abstract version presented in Theorem~\ref{Borel-Pompeiu_abstract}, we can conclude that the boundary operators $\mathfrak{B}_{1}(f,g,\mathbb{Z}^{7},h,0)$ and $\mathfrak{B}_{1}(f,g,\mathbb{Z}^{7},h,0)$ are, in fact, the discrete $F$-operators (or Cauchy transform) for octonions and split-octonions, respectively. Furthermore, the abstract Stokes' operator after plugging-in the discrete fundamental solution reduces to the discreet $T$-operator and to a function $f$ by using the properties of $E_{h}^{+}$. Hence, the unified approach is particularly useful after defining basic operators explicitly, because all function-theoretic formulae can be defined concisely. But, nonetheless, the use of discrete operator calculus in practice requires explicit formulae.\par
We finish this section by presenting the discrete octonionic and split-octonionic Cauchy formulae for the half-spaces:
\begin{theorem}
Let $f$ be a discrete left backward monogenic function wit respect to the operator $D_{h}^{-}$, and let $E_{h}^{+}$ be th discrete fundamental solution to the discrete Cauchy-Riemann operator $D_{h}^{+}$. Then the discrete octonionic Cauchy formula for the upper half-space $h\mathbb{Z}_{+}^{8}$ is given by
\begin{equation*}
\mathcal{A}_{1}(f,g,\mathbb{Z}_{+}^{8},I) - \mathfrak{B}_{1}(f,g,\mathbb{Z}^{7},h,0) = \left\{\begin{array}{cl}
0, & m\notin \mathbb{Z}_{+}^{8}, \\
f(mh), & m\in \mathbb{Z}_{+}^{8},
\end{array} \right.
\end{equation*}
and the discrete split-octonionic Cauchy formula for the upper half-space $h\mathbb{Z}_{+}^{8}$ is given by
\begin{equation*}
\mathcal{A}_{1}(f,g,\mathbb{Z}_{+}^{8},J) - \mathfrak{B}_{1}(f,g,\mathbb{Z}^{7},h,0) = \left\{\begin{array}{cl}
0, & m\notin \mathbb{Z}_{+}^{8}, \\
f(mh), & m\in \mathbb{Z}_{+}^{8}.
\end{array} \right.
\end{equation*}
These formulae hold for all discrete functions $f$ and $g$ under the condition that the series converge. The formulae for the lower half-space $h\mathbb{Z}_{-}^{8}$ are defined analogously, where the boundary layers $0$ and $-1$ need to be used. 
\end{theorem}\par
Finally we would like to remark, that discrete Hardy spaces can also easily be defined in the split-octonionic case. The procedure is completely analogously to the one described in \cite{CKKS,Cerejeiras} for the case of discrete Clifford analysis and in \cite{Krausshar_2} for the classical octonionic case. Therefore, this discussion is omitted here to avoid unnecessary duplications.\par

%
%

\end{document}